\documentstyle{amsppt}
\magnification=\magstep1
\baselineskip=14pt
\parindent=0pt
\parskip=14pt


\vsize=7.7in
\voffset=-.4in
\hsize=5.5in



\overfullrule=0pt

\define\A{{\Bbb A}}
\define\C{{\Bbb C}}

\define\R{{\Bbb R}}
\define\Q{{\Bbb Q}}

\define\Z{{\Bbb Z}}


\redefine\H{\frak H}




\redefine\b{\beta}

\redefine\l{\lambda}
\redefine\o{\omega}
\define\ph{\varphi}

\define\s{\sigma}
\redefine\P{\Phi}
\predefine\Sec{\S}

\redefine\S{{\Cal S}}



\define\back{\backslash}

\define\lra{\longrightarrow}
\redefine\tt{\otimes}
\define\scr{\scriptstyle} 
\define\liminv#1{\underset{\underset{#1}\to\leftarrow}\to\lim}
\define\limdir#1{\underset{\underset{#1}\to\rightarrow}\to\lim}

\define\isoarrow{\ {\overset{\sim}\to{\longrightarrow}}\ }

\define\nass{\noalign{\smallskip}}


\define\CH{\widehat{\roman{CH}}}
\define\ZH{\widehat{\Cal Z}}  


\redefine\O{\Omega}
\predefine\oldvol{\vol}
\redefine\vol{\text{\rm vol}}
\define\pr{\text{\rm pr}}


\define\bom{\o}


\redefine\div{\text{\rm div}}


\define\Spec{\text{\rm Spec}\,}

\define\Sym{\text{\rm Sym}}
\define\tr{\text{\rm tr}}

\font\cute=cmitt10 at 12pt
\font\smallcute=cmitt10 at 9pt
\define\kay{{\text{\cute k}}}
\define\smallkay{{\text{\smallcute k}}}

\define\sig{\text{\rm sig}}

%

%

\redefine\Re{\text{\rm Re}}
\define\SL{\text{\rm SL}}

\define\und#1{\underline{#1}}



\define\CT#1{\operatornamewithlimits{CT}_{#1}}

\define\Pic{\text{\rm Pic}}

\define\Pich{\widehat{\Pic}}
\define\degh{\widehat{\deg}\ }

\define\tent#1{ \vphantom{\vbox to #1pt{}} }   

\define\undd#1{\und{\und{#1}}}


\font\loud=cmb10 at 14pt   

\font\quote=lcmssi8


\define\boecherer{\bf1}
\define\bocherersptriple{\bf2}
\define\borchinventI{\bf3}
\define\borchinventII{\bf4}
\define\borcherdsK{\bf5}  
\define\borchduke{\bf6}
\define\borchdukeII{\bf7}
\define\bostbourb{\bf8}
\define\bost{\bf9}
\define\bostumd{\bf10}
\define\bostgilletsoule{\bf11}
\define\bruinierI{\bf12}
\define\bruinierII{\bf13}
\define\brkuehn{\bf14}
\define\bkk{\bf15}  
\define\cohen{\bf16}
\define\deligne{\bf17}
\define\freitaghermann{\bf18}
\define\funkethesis{\bf19}
\define\funkecompo{\bf20}
\define\funkemillson{\bf21}
\define\garrettII{\bf22}
\define\garrett{\bf23}
\define\gsihes{\bf24}
\define\GN{\bf25}
\define\grossquasi{\bf26}
\define\grossmot{\bf27}
\define\grosskeating{\bf28}
\define\grosskohnenzagier{\bf29}
\define\grosskudla{\bf30}
\define\grosszagier{\bf31}
\define\harris{\bf32}
\define\harriskudlaII{\bf33}
\define\harriskudlaIII{\bf34}
\define\harrisksweet{\bf35}
\define\hirzebruchzagier{\bf36}
\define\howeps{\bf37}
\define\kitaoka{\bf38}
\define\kitaokatriple{\bf39}
\define\hirschberg{\bf40}
\define\duke{\bf41}
\define\annals{\bf42}
\define\bourbaki{\bf43}
\define\Bints{\bf44}
\define\icm{\bf45}
\define\CDM{\bf46}
\define\kmI{\bf47}
\define\kmII{\bf48}
\define\kmihes{\bf49}
\define\kmcana{\bf50}
\define\krannals{\bf51}
\define\krHB{\bf52}
\define\krinvent{\bf53}
\define\krsiegel{\bf54}
\define\tiny{\bf55}  
\define\kryII{\bf56}   
\define\kryIII{\bf57}  
\define\kuehn{\bf58}
\define\li{\bf59}
\define\loojenga{\bf60}
\define\maillotroessler{\bf61}
\define\mcgraw{\bf62}
\define\milne{\bf63}
\define\niwa{\bf64}
\define\psrallis{\bf65}
\define\psrallistriple{\bf66}
\define\dprasad{\bf67}
\define\rallisinnerprod{\bf68}
\define\satake{\bf69}
\define\shimurahalf{\bf70}
\define\shimuraorth{\bf71}
\define\shimuraconf{\bf72}
\define\shintani{\bf73}
\define\soulebook{\bf74}
\define\vdgeer{\bf75}   
\define\vdgeerbook{\bf76}
\define\waldspurger{\bf77}
\define\weilactaII{\bf78}
\define\yangden{\bf79}
\define\yangiccm{\bf80}
\define\yangMSRI{\bf81}
\define\zagierII{\bf82}
\define\zagier{\bf83}   

\define\hfb{\hfill\break}

\define\hbom{\widehat{\bom}}

\define\wphi{\widehat{\phi}}
\define\Pet{{\roman{Pet}}}
\define\thh{\hat\theta}

\define\qeq{\ \overset{??}\to{=}\ }

\font\footbold=msbm10 at 8pt    
\define\Bbbf#1{\text{\footbold#1}}

\centerline{\bf Special cycles and derivatives of Eisenstein series}
\medskip
\centerline{Stephen S. Kudla\footnote{Partially 
supported by NSF grant DMS-9970506 and by a Max-Planck Research Prize 
from the Max-Planck Society and Alexander von Humboldt Stiftung. }}
\medskip
\medskip 
\centerline{\hfill  
$\matrix\text{\quote``A man hears what he wants to hear}\\\text{\quote and disregards the
rest."}\\
\scr\text{Simon and Garfunkel}\\
\scr\text{\it The boxer}\endmatrix$}

\medskip
\medskip
This article is an expanded version of a lecture given at the conference on 
Special Values of Rankin L-Series at MSRI in December of 2001. I have tried to retain 
some of the tone of an informal lecture. In particular, I have attempted to 
outline, in very broad terms, a program involving relations among:
\roster
\item"{(i)}" algebraic cycles
\item"{(ii)}" Eisenstein series and their derivatives
\item"{(iii)}" special values of Rankin-Selberg L-functions and their derivatives,
\endroster
ignoring many important details and serious technical problems in the process. 
I apologize at the outset for the very speculative nature of the picture given here. 
I hope that, in spite of many imprecisions, the sketch will provide a context for a variety 
of particular cases where precise results 
have been obtained. Recent results on one of these, part of an ongoing 
joint project with Michael Rapoport and Tonghai Yang on which much of the conjectural 
picture is based, are described in Yang's article in this volume, \cite{\yangMSRI}. A less speculative 
discussion of some of this material can be found in \cite{\bourbaki}, \cite{\icm}, and \cite{\CDM}. 

I would like to thank my collaborators B. Gross, M. Harris, J. Millson, S. Rallis, 
M. Rapoport and  T. Yang for generously sharing their mathematical ideas and for 
their support over many years.  I would also like to thank R. Borcherds, J.-B. Bost, 
J. Cogdell, J. Funke, R. Howe, D. Kazhdan, K. Keating, J. Kramer, U. K\"uhn,  J.S. Li, J.
Schwermer, and D. Zagier for helpful discussions, comments and suggestions.  Finally, I
would like to thank Henri Darmon and Shou-Wu Zhang for organizing  such an enjoyable and
inspiring program and for their very helpful suggestions which
significantly improved this article. Finally, I would like to thank MSRI for its ever excellent hospitality. 
\vskip .5in
\vfill
\eject

\centerline{\loud I}
\medskip
\centerline{\loud An attractive family of varieties}

\subheading{\Sec1. Shimura varieties of orthogonal type}

We begin with the following data:\footnote{Recall that, if $C(V) = C^+(V)\oplus C^-(V)$ is the Clifford algebra 
of $V$ with its $2$--grading, then there is a canonical embedding $V\hookrightarrow C^-(V)$, and
$$\roman{GSpin}(V)= \{\ g\in C^+(V)^\times \mid gg^\iota = \nu(g),\ \text{and}\ gVg^{-1} = V\ \},$$
where $\iota$ is the main involution of $C(V)$ and $\nu(g)$ is a scalar. There is an exact sequence
$$1\lra Z \lra G\lra SO(V)\lra 1$$
and the spinor norm homomorphism $\nu:G\rightarrow \Bbbf G_m$. For $n\ge 1$, strong approximation holds for
the  semi-simple, simply connected group $\roman{Spin}(V)
=\ker(\nu)$. } 
$$\align
V, \ (\ ,\ )&=\text{inner product space over $\Q$}\\
\nass
\sig(V)&=(n,2)\\
\nass
G&=\roman{GSpin}(V)\tag1.1\\
\nass
D&
= \{\ w\in V(\C)\mid (w,w)=0,\ (w,\bar w)<0\ \}/\C^\times  \ \subset\  \Bbb P(V(\C))\\
\nass
n&=\dim_{\,\C} D.\\
\endalign
$$
This data determines a Shimura variety $M = \roman{Sh}(G,D)$, with a canonical model over 
$\Q$, where, for $K\subset G(\A_f)$ a compact open subgroup, 
$$M_K(\C)\simeq G(\Q)\back \bigg( D\times G(\A_f)/K\bigg).\tag1.2$$
Note that $D = D^+\cup D^-$ is a union of two copies of a bounded domain of type IV, \cite{\satake}, p.285.
They are interchanged by the complex conjugation $w\mapsto \bar w$.  
If we let $G(\R)^+$ be the subgroup of $G(\R)$ which preserves $D^+$ and write
$$G(\A) = \coprod_j G(\Q)G(\R)^+ g_jK,\tag1.3$$
then 
$$M_K(\C) \simeq \coprod_j \Gamma_j\back D^+,\tag1.4$$
where $\Gamma_j = G(\Q)\cap G(\R)^+g_jKg_j^{-1}$. 
Thus, for general $K$, the quasi-projective variety $M_K$ can have many components\footnote{In fact,
for $n\ge1$, 
by strong approximation, 
$$\pi_0(M_K) \simeq \Bbbf Q^\times\Bbbf R^\times_+\back \Bbbf A^\times/\nu(K),$$
where $\nu(K)$ is the image of $K$ under the norm map $\nu$.
}
and  
the individual components are only rational over some cyclotomic extension. The action 
of the Galois group on the components is described, for example, in \cite{\deligne}, 
\cite{\milne}. 

$M_K$ is quasi-projective of dimension $n$ over $\Q$, and projective if and only 
if the rational quadratic space $V$ is anisotropic. By Meyer's Theorem, this can 
only happen for $n\le2$. 
In the range $3\le n\le 5$, we can have $\roman{witt}(V)=1$, where $\roman{witt(V)}$ 
is the dimension of a maximal isotropic 
$\Q$-subspace of $V$. For $n\ge6$, $\roman{witt(V)}=2$. 
A nice description of the Baily--Borel compactification of $\Gamma\back D^+$ 
and its
toroidal desingularizations can be found in \cite{\loojenga}. 

For small values of $n$, the $M_K$'s include many classical varieties, for example:
\roster
\item"{}" $n=1$, modular curves and Shimura curves, \cite{\annals},
\item"{}" $n=2$, Hilbert-Blumenthal surfaces and quaternionic versions, \cite{\krHB}, \cite{\vdgeerbook},
\item"{}" $n=3$, Siegel 3-folds and quaternionic analogues, \cite{\krsiegel}, \cite{\vdgeer}, \cite{\GN}, 
\item"{}" $n\le 19$, moduli spaces of K3 surfaces,  \cite{\borcherdsK}. 
\endroster
Of course, such relations are discussed in many places, cf., for example, 
\cite{\freitaghermann}. 

The most familiar example arises for the three dimensional rational quadratic space
$$V=\{\ x\in M_2(\Q)\mid \tr(x)=0\ \},$$
with $Q(x) = \det(x)$. 
In this case, $G=\roman{GSpin}(V) =\roman{GL}(2)$, the 
homomorphism $G\rightarrow SO(V) \simeq \roman{PGL}(2)$ is defined via the conjugation 
action on $V$, $D\simeq \H^+\cup \H^-$, the union of the upper and lower half planes, 
and the $M_K$'s are the usual modular curves. 
Note that it is more
convenient to work with
$\roman{GL}(2)$  rather than $\roman{PGL}(2)$, and this leads to the choice of 
$\roman{GSpin}(V)$ rather than $SO(V)$ or the disconnected group $O(V)$ in general. 
More details about the case of signature $(1,2)$ are given in the Appendix below.

More generally, one could consider quadratic spaces $V$ over 
a totally real field $\kay$ 
with $\sig(V_{\infty_i}) = (n,2)$ for $\infty_i\in S_1$ and $\sig(V_{\infty_i})=(n+2,0)$ 
for $\infty_i\in S_2$ where $S_1\cup S_2$ is a disjoint decomposition 
of the set of archimedean places of $\kay$. 
If $S_2\ne\emptyset$, then the varieties $M_K$ are always projective. 
Such compact quotients are considered in \cite{\kmI}, \cite{\kmII}, and \cite{\duke}. 
For a
discussion of  automorphic forms in this situation from a classical point of view, see
\cite{\shimuraorth}.  
\medskip

\subheading{\Sec2. Algebraic cycles}

An attractive feature of this family of Shimura varieties is that they have 
many algebraic cycles; in fact, there are sub-Shimura varieties of the same type 
of all codimensions. These can be 
constructed as follows.

Let $\Cal L_D$ be the homogeneous line bundle over $D$ with
$$\Cal L_D \setminus \{0\} = \{w\in V(\C)\mid (w,w)=0, \ (w,\bar w)<0\},\tag2.1$$
so that $\Cal L_D$ is the restriction to $D$ of the bundle $\Cal O(-1)$ 
on $\Bbb P(V(\C))$. 
We equip $\Cal L_D$ with the hermitian metric $||\ ||$ given by $||w||^2 = |(w,\bar w)|$.
The action of $G(\R)$ on $D$ lifts in a natural way to an action on $\Cal L_D$, 
and hence, this bundle descends to a line bundle $\Cal L$ on 
the Shimura variety $M$. For example, 
for a given compact open subgroup $K$, $\Cal L_K\rightarrow M_K$, has 
a canonical model over $\Q$, \cite{\harris}, \cite{\milne}, 
and
$$\Cal L_K(\C) \simeq G(\Q)\back \bigg( \Cal L_D\times G(\A_f)/K\bigg).\tag2.2$$

Any rational vector $x\in V(\Q)$ defines a section $s_x$ over $D$ of 
the dual bundle $\Cal L_D^\vee$
by the formula
$$(s_x,w) = (x,w),\tag2.3$$
and, for $x\ne0$, the (possibly empty) divisor\footnote{a rational quadratic divisor 
in Borcherds' terminology, \cite{\borchinventI}.} in $D$ of this section is given by
$$\div(s_x) = \{w\in D\mid (x,w)=0\}/\C^\times =: D_x \subset D.\tag2.4$$
Assuming that $Q(x):=\frac12(x,x)>0$ and setting
$$\align 
V_x &= x^\perp\tag2.5\\
\noalign{and}
G_x&=\roman{GSpin}(V_x) = \text{stabilizer of $x$ in $G$,}\\
\noalign{there is a sub-Shimura variety}
Z(x)&: \roman{Sh}(G_x,D_x) \lra \roman{Sh}(G,D)= M\tag2.6
\endalign
$$
giving a divisor $Z(x)_K$, rational over $\Q$, on $M_K$ for each $K$. 

If $Q(x)\le 0$, and $x\ne 0$, then the section $s_x$ is never zero on $D$,  
so that $D_x=\emptyset$.
If $x=0$, then we formally set $D_x=D$ and take $Z(0) = M$.

More generally, given an $r$--tuple of vectors $x\in V(\Q)^r$ 
define
$V_x$, $G_x$, $D_x$ by the same formulas.  
If the span $\und{x}$ of the components of $x$ has dimension $r(x)$ and if the matrix
$$Q(x) = \frac12\big((x_i,x_j)\big)\tag2.7$$
is positive semidefinite of rank $r(x)$, then the restriction of $(\ ,\ )$ to $V_x$ 
has signature $(n-r(x),2)$, and  
there is a corresponding 
cycle $Z(x):Sh(G_x,D_x)\rightarrow Sh(G,D)=M$, of codimension
$r(x)=\roman{rk}(Q(x))\le r$. If the rank of $Q(x)$ is less than $r(x)$ or if $Q(x)$ 
is not positive semi-definite, then $Z(x)=\emptyset$, \cite{\kmI}, \cite{\kmII}, \cite{\kmihes}. 

For $g\in G(\A_f)$, we can also make a `translated' cycle $Z(x,g)$ 
where, at level $K$, 
$$\align
Z(x,g;K): G_x(\Q)\back \bigg( D_x\times G_x(\A_f)/ K_x^g\bigg)  &\lra
G(\Q)\back
\bigg( D\times G(\A_f)/K\bigg)=M_K(\C),\tag2.8\\
\nass
G_x(\Q)(z,h) K_x^g\quad&\mapsto\quad G(\Q)(z,hg)K.\endalign
$$
where we write $K_x^g = G_x(\A_f)\cap gKg^{-1}$ for short. 
This cycle is again rational over $\Q$.

Finally, we form certain weighted combinations of these cycles, essentially by summing over 
integral $x$'s with a fixed matrix of inner products, 
\cite{\duke}. More precisely, 
suppose that a $K$-invariant Schwartz function\footnote{For example, for $r=1$, $\ph$ 
might be
the characteristic function of the closure in $V(\Bbbf A_f)$ of a coset $\mu+L$ 
of a lattice $L\subset V$.} 
$\ph\in
S(V(\A_f)^r)^K$ on
$r$ copies of the finite adeles 
$V(\A_f)$ of $V$ and  $T\in \Sym_r(\Q)_{\ge0}$ are given. 
Let 
$$\O_T = \{\ x\in V^r\mid Q(x) = T\ \},\tag2.9$$
and, assuming that $\O_T(\Q)$ is nonempty, fix an element $x\in\O_T(\Q)$ and write
$$\O_T(\A_f) \cap \roman{supp}(\ph)= \coprod_j K g_j^{-1} x\tag2.10$$
for $g_j\in G(\A_f)$. Note that the sum is finite since the $K$-orbits give an open cover 
of the compact 
set $\O_T(\A_f)\cap \roman{supp}(\ph)$. Then
there is a cycle $Z(T,\ph;K)$ in $M_K$ defined by 
$$Z(T,\ph;K) = \sum_j \ph(g_j^{-1}x)\,Z(x,g_j;K)\tag2.11$$
of codimension $\roman{rank}(T) =:r(T)$, 
given by a weighted combination of the $Z(x)$'s for $x$ with $Q(x)=T$. 

These weighted 
cycles have nice properties, 
\cite{\duke}. For example, if $K'\subset K$ and\hfb $\roman{pr}:M_{K'}\rightarrow
M_K$ is the corresponding covering map, then
$$\roman{pr}^*Z(T,\ph;K) = Z(T,\ph;K').\tag2.12$$
Thus it is reasonable to drop $K$ from the notation and write simply $Z(T,\ph)$. 

{\bf Example.} The classical Heegner divisors, traced down to $\Q$, arise in the 
case $n=1$, $r=1$. A detailed description is given in Appendix I below.  

\subheading{\Sec3. Modular generating functions}

In this section, we discuss the generating functions which can be constructed 
from the cycles $Z(T,\ph)$, by taking their classes either in cohomology 
or in Chow groups. The main goal is to prove that such generating functions 
are, in fact, modular forms. Of course, these constructions are modeled
on the work of Hirzeburch and Zagier \cite{\hirzebruchzagier} on generating functions 
for the cohomology classes of curves on Hilbert-Blumenthal surfaces.

{\bf 3.1. Classes in cohomology.}\hfb
The cycles defined above are very special cases of the locally symmetric cycles 
in Riemannian locally symmetric spaces studied some time ago in a long collaboration with John
Millson, 
\cite{\kmI}, \cite{\kmII}, \cite{\kmihes}. The results described in this section 
are from that joint work. 
For $T\in \Sym_r(\Q)_{\ge0}$ and a weight function $\ph$, there are cohomology classes
$$[\,Z(T,\ph)\,]\in H^{2r(T)}(M_K)\qquad\text{and}
\qquad [\,Z(T,\ph)\,]\cup [\Cal L^\vee]^{r-r(T)}\in H^{2r}(M_K),\tag3.1$$
where $r(T)$ is the rank of $T$ and $[\Cal L^\vee]\in
H^2(M_K)$ is the cohomology  class of the dual $\Cal L^\vee$ of the line bundle $\Cal L$. 
Here we view our cycles as defining linear functionals on the 
space of compactly supported closed forms, and hence these classes lie in 
the absolute cohomology $H^\bullet(M_K)$ of $M_K(\C)$ with complex coefficients. 

In \cite{\kmihes}, we proved:
\proclaim{Theorem 3.1} 
For $\tau = u+iv\in \H_r$, the Siegel space of genus $r$, the holomorphic function
$$\phi_r(\tau,\ph) = \sum_{T\in \Sym_r(\Q)_{\ge0}}[\,Z(T,\ph)\,]\cup 
[\Cal L^\vee]^{r-r(T)}\,q^T,$$
is a Siegel modular form of genus $r$ and weight $\frac{n}2+1$ valued in $H^{2r}(M_K)$. \hfb
Here $q^T=e(\tr(T\tau))$. 
\endproclaim
\demo{Idea of Proof} The main step is to construct a theta 
function {\it valued in the closed $(r,r)$--forms on $M_K(\C)$}. Let $A^{(r,r)}(D)$ 
be the space of smooth $(r,r)$--forms on $D$, and let $S(V(\R)^r)$ 
be the Schwartz space of $V(\R)^r$. The group $G(\R)$ acts naturally on 
both of these spaces. For $\tau\in \H_r$, there is a 
Schwartz form,
\cite{\kmI}, \cite{\kmII}, \cite{\kmihes}, 
$$\ph_\infty^r(\tau)\in \bigg[\,S(V(\R)^r)\tt A^{(r,r)}(D)\,\bigg]^{G(\R)}\tag3.2$$
with the key property that, 
for all $x\in V(\R)^r$, $d\ph_\infty^r(\tau,x)=0$, i.e., $\ph_\infty^r(\tau,x)$
is  a closed form on $D$. 
Note that, for $g\in G(\R)$ and $x\in V(\R)^r$, the $G(\R)$--invariance in (3.2) means 
that
$$g^*\ph_\infty^r(\tau,x) = \ph_\infty^r(\tau,g^{-1}x).\tag3.3$$
Thus, for example, for a fixed $x$,  $\ph_\infty(\tau,x)\in A^{(r,r)}(D)^{G(\R)_x}$
is a closed $G(\R)_x$--invariant form on $D$. 
Note that $\ph_\infty^r(\tau)$ is {\it not} holomorphic in $\tau$. 
For any $\ph\in S(V(\A_f)^r)^K$,  
the Siegel theta function
$$\theta_r(\tau,\ph)(\cdot,h) := \sum_{x\in V(\Q)^r}\ph_\infty^r(\tau,x)\ph(h^{-1}x)\tag3.4$$ 
defines a closed $(r,r)$--form on $M_K(\C)$. By the 
standard argument based on Poisson summation, which is explained in a little more 
detail in Appendix II, $\theta_r(\tau,\ph)$ is modular of weight $\frac{n}2+1$
for a subgroup
$\Gamma'\subset
\roman{Sp}_r(\Z)$, depending on $\ph$.
Finally, the cohomology class 
$$\phi_r(\tau,\ph) = [\theta_r(\tau,\ph)]\tag3.5$$ 
of the theta form (3.4) coincides with the {\it holomorphic} generating function 
of the Theorem and hence this generating function is  
also modular of weight $\frac{n}2+1$. 
\qed\enddemo 

The Schwartz forms satisfy the cup product identity:
$$\ph_\infty^{r_1}(\tau_1) \wedge \ph_\infty^{r_2}(\tau_2) = 
\ph_\infty^{r_1+r_2}(\pmatrix \tau_1&{}\\{}&\tau_2\endpmatrix),\tag3.6$$
where the left side is an element of the space $S(V(\R)^{r_1})\tt S(V(\R)^{r_2})\tt 
A^{(r,r)}(D)$, with $r=r_1+r_2$, and $\tau_j\in \H_{r_j}$. Hence,
for weight functions $\ph_j\in S(V(\A_f)^{r_j})$,  
one has the
identity for the theta forms 
$$\theta_{r_1}(\tau_1,\ph_1)\wedge \theta_{r_2}(\tau_2,\ph_2) =
 \theta_{r}(\pmatrix \tau_1&{}\\{}&\tau_2\endpmatrix,\ph_1\tt\ph_2).\tag3.7$$
Passing to cohomology, (3.7) yields the pleasant identity, 
\cite{\duke}:
$$\phi_{r_1}(\tau_1,\ph_1)\cup \phi_{r_2}(\tau_2,\ph_2) = \phi_{r_1+r_2}
\big(\pmatrix \tau_1&{}\\{}&\tau_2\endpmatrix, \ph_1\tt\ph_2 \big),\tag3.8$$
for the cup product of the generating functions valued in $H^\bullet(M_K)$.
Comparing coefficients, we obtain the following formula for
the cup product of our classes.\hfb Suppose that $T_1\in \Sym_{r_1}(\Q)_{>0}$ 
and $T_2\in \Sym_{r_2}(\Q)_{>0}$. Then
$$[Z(T_1,\ph_1)]\cup [Z(T_1, \ph_2)] = \sum_{\matrix \scr T\in \Sym_r(\Q)_{\ge0}\\ 
\nass
\scr T = \pmatrix \scr T_1&\scr *\\ \scr*& \scr T_2\endpmatrix\endmatrix}
[Z(T,\ph_1\tt\ph_2)]\cup [\Cal L^\vee]^{r-\roman{rk}(T)}.\tag3.9$$

{\bf 3.2. Classes in Chow groups.}\hfb
We can also take classes of the cycles in the usual Chow groups\footnote{We only 
work with rational coefficients.}.
For this, when $V$ is anisotropic so that $M_K$ is compact, we consider
the  classes 
$$\{\,Z(T,\ph)\,\}\in \roman{CH}^{r(T)}(M_K)\qquad\text{and}
\qquad \{\,Z(T,\ph)\,\}\cdot \{\Cal L^\vee\}^{r-r(T)}\in \roman{CH}^{r}(M_K)\tag3.10$$
in the Chow groups of $M_K$, and corresponding generating functions
$$\phi^{\roman{CH}}_r(\tau,\ph)= \sum_{T\in \Sym_r(\Q)_{\ge0}}
\{\,Z(T,\ph)\,\}\cdot \{\Cal L^\vee\}^{r-r(T)}\,q^T\tag3.11$$
valued in $\roman{CH}^r(M_K)_\C$. Here $\cdot$ denotes the 
product in the Chow ring $\roman{CH}^{\bullet}(M_K)$, 
and $\{\Cal L^\vee\}\in \roman{Pic}(M_K)\simeq\roman{CH}^1(M_K)$ 
is the class of $\Cal L^\vee$. Note that, for the cycle
class map: 
$$cl: \roman{CH}^r(M_K)\lra H^{2r}(M_K),\tag3.12$$
we have
$$\phi^{\roman{CH}}_r(\tau,\ph)\mapsto \phi_r(\tau,\ph),\tag3.13$$
so that the generating function $\phi^{\roman{CH}}_r(\tau,\ph)$ `lifts'
the cohomology valued function 
$\phi_r(\tau,\ph)$, which is modular by Theorem~3.1. 

If $V$ is isotropic, let $\widetilde{M}_K$ be a smooth toroidal compactification 
of $M_K$, \cite{\loojenga}. Let $Y_K = \widetilde{M}_K \setminus M_K$ 
be the compactifying divisor, and let 
$\roman{CH}^1(\widetilde{M}_K,Y_K)$ be the quotient of $\roman{CH}^1(\widetilde{M}_K)$ 
by the subspace generated by the irreducible components of $Y_K$.  We use the 
same notation for the classes $\{\,Z(T,\ph)\,\}$ of our cycles in this group. 
\hfb
The following result is due to Borcherds, \cite{\borchduke}, \cite{\borchdukeII}.

\proclaim{Theorem 3.2}
For $r=1$ and for a $K$-invariant weight function 
$\ph\in S(V(\A_f))^K$,  the generating
function
$$\phi^{\roman{CH}}_1(\tau,\ph) = \{\Cal L^\vee\}\,\ph(0) + \sum_{t>0} 
\{\,Z(t,\ph)\,\}\, q^t$$ 
is an elliptic modular form of weight $\frac{n}2+1$ valued in
$\roman{CH}^1(\widetilde{M}_K,Y_K)$. 
\endproclaim

\demo{Proof} Since the result as stated is not quite in \cite{\borchduke}, 
we indicate the precise relation to Borcherds' formulation. For a lattice $L\subset V$
on which the quadratic form $Q(x)=\frac12(x,x)$ is integer valued, let 
$L^\vee =\{x\in V\mid (x,L)\subset \Z\}$ be the dual lattice. Let $S_L\subset
S(V(\A_f))$ be the finite dimensional subspace spanned by the characteristic 
functions $\ph_\l$ of the closures in $V(\A_f)$ of the cosets $\l+L$ where $\l\in
L^\vee$. Every $\ph\in S(V(\A_f))$ lies in some $S_L$ for sufficiently small $L$.
There is a (finite Weil) representation $\rho_L$ of a central extension $\Gamma'$ 
of $\SL_2(\Z)$ on $S_L$. Suppose that $F$ is a holomorphic function on $\H$, 
valued in $S_L$, which is modular of weight $1-\frac{n}2$, i.e., for all $\gamma'\in
\Gamma'$
$$F(\gamma'(\tau)) = j(\gamma',\tau)^{2-n}\rho_L(\gamma')F(\tau),\tag3.14$$
where $j(\gamma',\tau)$ with $j(\gamma',\tau)^2 = (c\tau+d)$ is the 
automorphy factor attached to $\gamma'$, and $\pmatrix
a&b\\c&d\endpmatrix$  is the projection of $\gamma'$ to $\SL_2(\Z)$. The function $F$ is
allowed to have  a pole of finite order at $\infty$, i.e., $F$ 
has a Fourier expansion of the form 
$$F(\tau) = \sum_{\l\in L^\vee/L}\sum_{m\in \Q} c_\l(m) \,q^m\,\ph_\l,\tag3.15$$
where only finitely many coefficients $c_\l(m)$ for $m<0$ can be nonzero. 
Note that, by the transformation law, $c_\l(m)$ can only be nonzero when 
$m\equiv -Q(\l)\mod \Z$. For any such $F$ where, in addition, all $c_\l(-m)$ 
for $m\ge 0$ are in $\Z$, Borcherd constructs a meromorphic function $\Psi(F)$
on $D$ with the properties, \cite{\borchduke}, \cite{\borchdukeII}, \cite{\borchinventI}:
\roster
\item"{(i)}" There is an integer $N$ such that for any $F$, $\Psi(F)^N$ is a 
meromorphic automorphic form of weight $k=N\,c_0(0)/2$, i.e., 
a meromorphic section of $\Cal L^{\tt k}$, and 
\item"{(ii)}" 
$$\div(\Psi(F)^2) = \sum_{\l\in L^\vee/L}\sum_{m>0} c_\l(-m)\,Z(m,\ph_\l).\tag3.16$$
\endroster
In \cite{\borchduke}, Borcherds defines a rational vector space $\roman{CHeeg}(M_K)$ with 
generators $y_{m,\l}$, for $\l\in L^\vee/L$ and $m>0$ with $m\equiv Q(\l)\mod \Z$, 
and $y_{0,0}$ and relations  
$$c_0(0)y_{0,0} + \sum_\l\sum_{m>0} c_\l(-m)\,y_{m,\l},\tag3.17$$
as $F$ runs over the quasi-modular forms of weight $1-\frac{n}2$, as above. 
Under the assumption that a certain space of vector valued forms has a basis with 
rational Fourier coefficients, Borcherds proved that 
the space $\roman{CHeeg}(M_K)$
is finite dimensional and that the generating function
$$\phi_1^B(\tau,L) = y_{0,0}\,\ph^\vee_0 + \sum_{\l}\sum_{m>0} y_{m,\l}\,q^m\,\ph_\l^\vee,\tag3.18$$
valued in $\roman{CHeeg}(M_K)\tt S_L^\vee$, is a modular 
form of weight $\frac{n}2+1$ for $\Gamma'$. Here $S_L^\vee$ is the dual space of 
$S_L$, with the dual representation $\rho_L^\vee$ of $\Gamma'$.  William McGraw \cite{\mcgraw}
recently proved that the necessary basis exists. 

To finish the proof of our statement, we choose a nonzero (meromorphic) section
$\Psi_0$ of $\Cal L$ and define a map
$$\align
\roman{CHeeg}(M_K) &\lra \roman{CH}(\widetilde{M}_K,Y_K),\\
\nass
y_{m,\l}&\mapsto Z(m,\ph_\l)\tag3.19\\
\nass
y_{0,0}&\mapsto -\div(\Psi_0). 
\endalign
$$
This is well defined, since a relation is mapped to 
$$-c_0(0)\,\div(\Psi_0) + \div(\Psi(F)^2) = N^{-1}\div(\Psi(F)^{2N}\,\Psi_0^{-2k}) \equiv
0,\tag3.20$$
since $\Psi(F)^{2N}\,\Psi_0^{-2k}$ is a meromorphic {\it function}
\footnote{Of course, one needs to check that it extends to 
a meromorphic function on $\widetilde{M}_K$.} on $M_K$. 
Since the generating function $\phi_1(\tau,\ph)$ is a finite linear combination 
of components of Borcherds' generating function, it is modular for some 
suitable subgroup of $\Gamma'$, as claimed. Note that, $\{\Cal
L^\vee\}=\{-\div(\Psi_0)\}$. 
\qed\enddemo

{\bf Problem 1:} Is $\phi^{\roman{CH}}_r(\tau,\ph)$ a Siegel modular form for $r>1$?  

{\bf Problem 2:} Does a cup product formula like (3.8) still hold?

{\bf Problem 3:} Define classes $\{\,\widetilde{\Cal L^\vee}\ \}\in \roman{Pic}(\widetilde{M}_K)$, 
$$\{\,\widetilde{Z}(T,\ph)\,\}\in \roman{CH}^{r(T)}(\widetilde{M}_K)\qquad\text{and}
\qquad \{\,\widetilde{Z}(T,\ph)\,\}\cdot \{\widetilde{\Cal L^\vee}\}^{r-r(T)}\in
\roman{CH}^{r}(\widetilde{M}_K)$$ so that the resulting generating function
$\phi^{\roman{CH}}_r(\tau,\ph)$  is modular.

Additional information about the map 
$$\roman{CHeeg}(M_K)/\Q\,y_{0,0} \lra 
\roman{CH}(\widetilde{M}_K,Y_K)/\Q\{\widetilde{\Cal L^\vee}\},\tag3.21$$
e.g., concerning injectivity,   
was obtained by Bruinier, \cite{\bruinierI}, \cite{\bruinierII}. 

\medskip
\medskip

\subheading{\Sec4. Connections with values of Eisenstein series}

To obtain classical scalar valued modular forms, one can apply 
linear functionals to the modular generating functions 
valued in cohomology. For a moment, we again assume that we are
in the case of compact quotient. Then, using the class
$$[\Cal L^\vee] \in H^2(M_K),\tag4.1$$
and the composition 
$$H^{2r}(M_K)\times H^{2(n-r)}(M_K) \lra H^{2n}(M_K) \overset{\deg}\to{\lra}\ \C,\tag4.2$$
of the cup product and the degree map, we have:
$$\deg(\,\phi_r(\tau,\ph)\cup [\Cal L^\vee]^{n-r}\,) =  \int_{M_K} 
\theta_r(\tau,\ph)\wedge \O^{n-r}:=I_r(\tau,\ph),\tag4.3$$
where $\O$ is the Chern form of the line bundle $\Cal L^\vee$ for its natural metric. 

Now, the Siegel--Weil formula, \cite{\weilactaII}, \cite{\krannals}, 
relates the integral $I_r(\tau,\ph)$
of a theta function determined by a Schwartz function $\ph\in S(V(\A_f)^r)^K$ to 
a special value of a Siegel Eisenstein series $E_r(\tau,s,\ph)$, also associated to
$\ph$.  The parameter $s$ in this Eisenstein series is normalized as in Langlands, 
so that there is a functional equation with respect to $s\mapsto -s$, and the 
halfplane of absolute convergence is $\Re(s)>\frac{r+1}2$. Note
that, to apply the Siegel--Weil formula, we must first  relate the integral of the theta
{\it form} occurring in (4.3)  to the ad\`elic integral of the theta {\it function}
occurring in the  Siegel--Weil theory, \cite{\Bints}, section 4. Hence, we obtain:

\proclaim{The volume formula}
In the case of compact quotient, \cite{\duke}, 
$$\deg(\,\phi_r(\tau)\cup [\Cal L^\vee]^{n-r}\,) 
\overset{\raise 5pt\hbox{$\scr(1)$}}\to{=}\ I_r(\tau,\ph) 
\overset{\raise 5pt\hbox{$\scr(2)$}}\to{=}\ \vol(M_K,\O^n)
\cdot E_r(\tau,s_0,\ph),\tag4.4$$ 
where
$$ s_0  = \frac{n+1-r}2.\tag4.5$$
\endproclaim

In fact, this formula should hold in much greater 
generality, i.e., when $V$ is isotropic. First of all, the theta integral 
is termwise convergent whenever Weil's condition $r<n+1-\roman{witt}(V)$ holds, and so
the identity (2) in (4.4) is then valid. 
The result of \cite{\kmcana} can be applied and the argument given in \cite{\Bints} for
the case
$r=1$ carries over to prove the following. 
\proclaim{Theorem 4.1} When
$r<n+1-\roman{witt}(V)$, there is an identity
$$\sum_{T\ge0} \vol(Z(T,\ph),\O^{n-r(T)})\,q^T  = \vol(M,\O^n)\cdot
E_r(\tau,s_0,\ph).$$ 
\endproclaim
It remains to give a cohomological interpretation of the left 
side of this identity in the noncompact case. 
 
Some sort of regularization of 
the theta integral, say by the 
method of \cite{\krannals}, is needed to obtain an extension 
of (2) to the range $r\ge n+1-\roman{witt}(V)$, i.e., 
to the cases $r=n-1$ and $n$ when $\roman{witt}(V)=2$ or the case $r=n$, 
if $\roman{witt}(V)=1$.
For example, in the case of modular curves, where $n=r=1$, it was 
shown by Funke \cite{\funkecompo} that 
the theta integral coincides with Zagier's nonholomorphic 
Eisenstein series of weight $\frac32$, \cite{\zagierII}. 
In this case, there are definitely (non-holomorphic!) correction terms which 
do not have an evident cohomological meaning, although they are 
consistent with a suitable 
arithmetic Chow group formulation, cf. Yang's article, \cite{\yangMSRI}. Recent work 
of Funke and Millson \cite{\funkemillson}
considered the pairing of the theta form with closed forms not of 
compact support in the case of arithmetic 
quotients of hyperbolic $n$-space. 

{\bf Examples:}

{\bf1.} If $n=1$ and $V$ is anistropic, so that $M=M_K$ is a Shimura curve over $\Q$, then
$$\vol(M_K)\cdot E_1(\tau,\frac12,\ph) = \deg(\phi_1(\tau,\ph)) = 
\vol(M_K,\O) + \sum_{t>0}\deg(Z(t,\ph))\,q^t\tag4.6$$
is a special value at $s=\frac12$ of an Eisenstein series of 
weight $\frac32$, and the $Z(t,\ph)$'s are 
Heegner type $0$--cycles on $M_K$, cf. Appendix I. This identity is 
described in more detail in \cite{\kryII}.

{\bf2.} If $n=2$ and $V$  has $\roman{witt}(V)=1$ or is anisotropic, 
so that $M_K$ is a Hilbert--Blumenthal surface for some real
quadratic field or a compact analogue,  then
$$\vol(M_K)\cdot E_1(\tau,1,\ph) = \deg(\phi_1(\tau,\ph)\cup \O) = 
\vol(M_K) + \sum_{t>0}\vol(Z(t,\ph),\O)\,q^t\tag4.7$$
is the special value at $s=1$ of an Eisenstein series of weight $2$, and the
$Z(t,\ph)$'s are Hirzebruch-Zagier type curves \cite{\vdgeer} on $M_K$.

{\bf3.} If $n=2$ and $V$ is anisotropic, then
$$\multline
\vol(M_K)\cdot E_2(\tau,\frac12,\ph) = \deg(\phi_2(\tau,\ph))\\
\nass
\nass
 = \vol(M_K) 
+\sum_{\scr T\in \Sym_2(\Q)_{\ge0}\atop\scr r(T)=1} 
\vol(Z(T,\ph),\O)\,q^T + \sum_{T>0}\deg(Z(T,\ph))\,q^T\qquad
\endmultline\tag4.8$$
is the special value at $s=\frac12$ of an Eisenstein 
series of weight $2$ and genus $2$, 
and, for $T>0$, the $Z(T,\ph)$'s are $0$--cycles.  
Gross and Keating, \cite{\grosskeating}, observed such a phenomenon 
in the split case as well. 

{\bf4.} If $n=3$, and for $V$ with $\roman{witt}(V)=2$, $M_K$ is a Siegel 
modular 3-fold. Then, for $r=1$, the $Z(t,\ph)$'s are combinations 
of Humbert surfaces, and the identity of Theorem~4.1 asserts that 
their volumes are the Fourier coefficients of an Eisenstein 
series of weight $\frac52$,  \cite{\vdgeerbook}, \cite{\Bints}. 

\vskip 1in
\centerline{\loud II}
\medskip
\centerline{\loud Speculations on the arithmetic theory}
\medskip 

The main idea is that many of the phenomena described above have an analogue 
in arithmetic geometry, where
the varieties $M$ are replaced by integral models $\Cal M$ over $\Spec(\Z)$, 
the cycles $Z(T,\ph)$ are replaced by arithmetic cycles on $\Cal M$, and the classes 
of these cycles are taken in 
arithmetic Chow groups $\CH^r(\Cal M)$, \cite{\gsihes}, \cite{\soulebook}. One could then
define a function $\wphi_r$ valued in $\CH^r(\Cal M)$, 
lifting the modular generating function $\phi_r$ valued in cohomology. 
The main goal would be  to prove the modularity of $\wphi_r$ and to find analogues of
the identities  discussed above,  where the values of the Eisenstein series
occurring in section 4 are  replaced by their derivatives, i.e., by the second terms in
their Laurent expansions. 

At this point, I am going to give an idealized picture which ignores many serious 
technical problems involving: (i) the existence of good integral models, (ii) bad reduction
and the possible bad behavior of cycles at such places, (iii) noncompactness, 
boundary contributions, 
(iv) extensions of the Gillet-Soul\'e theory \cite{\gsihes} of arithmetic 
Chow groups $\CH^r(\Cal M)$ 
to allow singular metrics, (cf, Bost, \cite{\bost}, K\"uhn, \cite{\kuehn}, 
Burgos-Kramer-K\"uhn, \cite{\bkk}),
(iv) suitable definitions of Green functions, etc., etc. \footnote{The fastidious reader may
want to stop here.}

Nevertheless, the idealized picture can serve as a guide and, 
with sufficient effort, one can obtain 
rigorous results in various particular cases, 
\cite{\annals}, \cite{\krHB}, \cite{\krinvent}, \cite{\krsiegel},  \cite{\tiny},
\cite{\kryII}.  In all of these
cases, we only consider a  good maximal compact subgroup $K$ and a specific weight function
$\ph$  determined by a nice lattice, so, in the discussion to follow, we will suppress both
$K$  and $\ph$ from the notation.

\subheading{\Sec5. Integral models and cycles}

Suppose that we have: 
$$\align
\Cal M&=\text{ a regular model of $M$ over $\Spec(\Z)$,}\\
\nass
\CH^{\bullet}(\Cal M)&= \text{its (extended) arithmetic Chow groups,}\\
\nass
\hbom&= \text{ extension of the metrized line bundle $\Cal L^\vee$ to $\Cal M$,}
\tag5.1\\
\nass
\hbom&\in \Pich(\Cal M) \simeq \CH^1(\Cal M)\\
\nass
\Cal Z(T)&=\text{ an extension of $Z(T)$ on $M$ to a cycle on $\Cal M$, so that}\\
\nass
\nass
&\matrix {}&\Cal Z(T)&\lra& \Cal M&{}&{}\\
\nass
{}&\uparrow&{}&\uparrow&{}&{}\\
\nass
Z(T)=&\Cal Z(T)_\Q&\lra&\Cal M_\Q&=M&\qquad\qquad\text{(generic fibers)}.\endmatrix
\endalign
$$

Finally, to obtain classes in the Gillet--Soul\'e arithmetic 
Chow groups $\CH^r(\Cal M)$ from the $\Cal Z(T)$'s we need 
Green forms, \cite{\gsihes}, \cite{\bostgilletsoule}, \cite{\bostbourb}. 
Based on the constructions for $r=1$, \cite{\annals}, 
and for $r=2$, $n=1$, \cite{\bourbaki}, we suppose that these have 
the following form:
$$\align
\tau &= u+iv\in \H_r\\
\nass
\Xi(T,v)&= \text{ Green form for $Z(T)$, depending on $v$,}\tag5.2\\
\nass
\ZH(T,v)& = (\Cal Z(T),\Xi(T,v))\in \CH^{r(T)}(\Cal M).
\endalign
$$
In all cases done so far, $0\le n\le 3$, \cite{\tiny},
 \cite{\annals}, \cite{\krinvent}, \cite{\krHB}, 
\cite{\krsiegel}, $M$ is of PEL type and the model $\Cal M$   
is obtained by extending the moduli problem over $\Q$ to a moduli 
problem over $\Spec(\Z)$ or, at least, $\Spec(\Z[N^{-1}])$ for a suitable 
$N$. The cycles $\Cal Z(T)$ are defined by imposing additional
endomorphisms satisfying various compatibilities, the special endomorphisms, 
cf. \cite{\bourbaki}, \cite{\icm} for further discussion. 

With such a definition, it can happen that $\Cal Z(T)$ is non-empty, even when 
the generic fiber $\Cal
Z(T)_\Q=Z(T)$ is empty.  For example, purely vertical divisors can occur in the 
fibers of bad reduction of the arithmetic surfaces attached to Shimura curves,
\cite{\krinvent}.  In addition, there can be cases where 
$\Cal Z(T)$ is empty, but $\Xi(T,v)$ is a nonzero smooth form on $M(\C)$, so that 
there are classes
$$\ZH(T,v) = (0,\Xi(T,v))\in \CH^r(\Cal M)\tag5.3$$
`purely vertical at infinity, even for $T$ not positive semi-definite, 
\cite{\annals}, \cite{\bourbaki}.

Finally, we define the {\bf arithmetic theta function:}
$$\widehat{\phi}_r(\tau) = \sum_{T\in \Sym_r(\Q)} \ZH(T,v)\cdot \hbom^{r-r(T)}\,
q^T \in \CH^r(\Cal M),\tag5.4$$
where $\cdot$ denotes the product in the arithmetic Chow ring $\CH^{\bullet}(\Cal M)$. 
Note that this function is not holomorphic in $\tau$, since the 
Green forms depend on $v$. Under the restriction
maps 
$$\roman{res}:\CH^r(\Cal M) \lra \roman{CH}^r(\Cal M_\Q) \lra H^{2r}(M),\tag5.5$$
we have 
$$\widehat{\phi}_r(\tau)\mapsto \phi^{\roman{CH}}_r(\tau)\mapsto \phi_r(\tau),\tag5.6$$
so that $\widehat{\phi}_r$ lifts $\phi^{\roman{CH}}_r$ and $\phi_r$ to the 
arithmetic Chow group. 

{\bf Problem 4:} Can the definitions be made so that $\widehat{\phi}_r(\tau)$ 
is a Siegel modular form of weight $\frac{n}2+1$ valued in 
$\CH^r(\Cal M)$, lifting $\phi_r$ and $\phi^{\roman{CH}}_r$?

At present, this seems out of reach, especially for $1<r<n+1$. 

{\bf Problem 5:} Is there an intersection product formula for the arithmetic Chow ring:
$$\widehat{\phi}_{r_1}(\tau_1)\cdot \widehat{\phi}_{r_2}(\tau_2) 
\overset{??}\to{=}\ \widehat{\phi}_{r_1+r_2} 
\big(\pmatrix \tau_1&{}\\{}&\tau_2\endpmatrix\big)\tag5.7$$
lifting the cup product relation (3.8) in cohomology?

\subheading{\Sec6. Connections with derivatives of Eisenstein series}

As in the standard Gillet--Soul\'e theory, suppose\footnote{Recall that, 
in the noncompact cases, we have to 
use some extended theory of arithmetic Chow groups, \cite{\bkk}, which allows the
singularities of the  natural metric on $\hbom$.} that
there is an arithmetic degree map
$$\degh:\CH^{n+1}(\Cal M) \lra \C,\tag6.1$$
and a height pairing
$$\langle\ ,\ \rangle:\CH^r(\Cal M)\times \CH^{n+1-r}(\Cal M) \lra \C, 
\qquad \langle \ZH_1,\ZH_2\rangle = \degh\!(\ZH_1\cdot\ZH_2).\tag6.2$$ 
These can be used to produce `numerical' generating functions from 
the $\widehat{\phi}_r$'s.

Let 
$$\Cal E_r(\tau,s) = C(s) \,E_r(\tau,s,\ph_0)\tag6.3$$
be the Siegel--Eisenstein series of weight $\frac{n}2+1$ and genus $r$ associated to
$\ph_0$, our standard weight function, with 
suitably normalizing factor $C(s)$, cf.~\cite{\kryII} for an example
of this normalization. The choice of $C(s)$ becomes important
in the cases in which  the leading term in nonzero. Then the following 
{\bf arithmetic volume formula} is an analogue  
of the volume formula of Theorem~4.1 above:

{\bf Problem 6:} For a suitable definition of $\Cal E_r(\tau,s)$, show that
$$\align
\Cal E'_r(\tau, s_0) &\qeq \langle\, \widehat{\phi}_r(\tau),\hbom^{n+1-r}\,\rangle\\
\nass
{}&=\ \sum_{T} \degh(\, \ZH(T,v)\cdot \hbom^{n+1-r(T)}\,)\,q^T.\tag6.4
\endalign
$$
where $
s_0=\frac{n+1-r}2$ is the 
critical value of $s$ occurring in the Siegel--Weil formula. 
Here $r$ lies in the range $1\le r\le n+1$.

{\bf Remarks:}\hfb
{\bf (i)} The identity (6.4) can be proved without knowing that $\wphi_r$ is modular, 
and one can obtain partial results by identifying corresponding Fourier coefficients 
on the two sides.\hfb 
{\bf (ii)} One can view the quantities $\degh\!(\, \ZH(T,v)\cdot \hbom^{n+1-r(T)}\,)$ 
as arithmetic volumes\hfb or heights, \cite{\bostgilletsoule}.  \hfb
{\bf (iii)} Assuming that $C(s_0) = \vol(M)$, the leading term 
$$\Cal E_r(\tau,s_0)= \vol(M)\,E_r(\tau,s_0)\tag6.5$$ 
of the normalized Eisenstein series at $s=s_0$ is just the generating function for geometric
volumes, via Theorem~4.1.\hfb {\bf(iv)} In the 
case $r=n+1$, so that $\wphi_r(\tau)\in \CH^{n+1}(\Cal M)$, the image of 
$\wphi_{n+1}$ in cohomology or in the usual Chow ring of $\Cal M_\Q$ 
is identically zero, since this group vanishes. On the other hand, 
the Eisenstein series $E_{n+1}(\tau,s)$ is incoherent in the sense of 
\cite{\annals}, \cite{\bourbaki}, the 
Siegel--Weil point is $s_0=0$, and $E_{n+1}(\tau,0)$ is also identically zero.
Thus the geometric volume identity is trivially valid. 
The arithmetic volume formula would then be
$$\degh\!(\wphi_{n+1}(\tau)) \qeq \Cal E'_{n+1}(\tau,0).\tag6.6$$

{\bf Examples:}

{\bf1. Moduli of CM elliptic curves.} $n=0$, $r=1$, \cite{\tiny}. Here $V$ is a negative
definite binary quadratic form given by the negative of the norm form of an imaginary 
quadratic field $\kay$, and $\Cal M$ is the moduli stack of elliptic curves
with CM by $O_{\smallkay}$, the ring of integers of $\kay$. For $t\in \Z_{>0}$, 
the cycle $\Cal Z(t)$ is either empty or is a $0$-cycle supported in 
a fiber $\Cal M_p$ for a prime $p$ determined by $t$. The
identity
$$\degh\!(\widehat{\phi}_1(\tau))= \Cal E_1'(\tau,0)\tag6.7$$
for the central derivative of an incoherent Eisenstein series of weight $1$ is 
proved in \cite{\tiny}, in the case in which $\kay$ has prime discriminant. 
The computation of the arithmetic degrees is based on the result 
of Gross, \cite{\grossquasi}, which is also the key to the geometric 
calculations in \cite{\grosszagier}. 

{\bf Remark:} In the initial work on the arithmetic situation \cite{\annals} 
and in the subsequent joint papers with Rapoport, \cite{\krsiegel} and 
\cite{\krHB}, the main idea was to view the central derivative of the 
incoherent Eisenstein series, restricted to the diagonal, as giving the 
height pairing of cycles in complementary 
degrees, cf. formula (6.10) below.  At the Durham conference in 1996, Gross
insisted that  it would be interesting to consider the `simplest case', i.e., $n=0$.
Following his suggestion, we obtained the results of \cite{\tiny} and 
came to see that the central derivative should {\it itself\,} have a nice geometric
interpretation, 
as a generating function for the arithmetic degrees of $0$--cycles on $\Cal M$,  
without restriction to the diagonal. This was a crucial step 
in the development of the picture discussed here.  

{\bf2. Curves on arithmetic surfaces.} 
$n=1$, $r=1$. Here, if $V$ is the
space of trace zero elements  of an indefinite division quaternion algebra over $\Q$, 
$\Cal M$ is the arithmetic surface associated to a Shimura curve.
For $t\in \Z_{>0}$, the cycle $\Cal Z(t)$ is a divisor on $\Cal M$ and can 
have vertical components. 
The identity
$$\langle\,\widehat{\phi}_1(\tau), \hbom\,\rangle=\Cal E_1'(\tau, \frac12)\tag6.8$$
is proved in \cite{\kryII}. Here $\Cal E_1(\tau,s)$ is a normalized Eisenstein 
series of weight $\frac32$. An unknown
constant, conjectured to be zero,  occurs in the definition of the class $\ZH(0,v)$ in
the  constant term of the generating function. 
This constant arises because we do not have, at present, an explicit 
formula for the quantity
$\langle\hbom,\hbom\rangle$ for the arithmetic surface attached to a Shimura curve. 
In the analogous example for modular curves, discussed in 
Yang's talk, \cite{\yangMSRI}, the quantity $\langle\hbom,\hbom\rangle$ is known, thanks
to the work of  Ulf K\"uhn \cite{\kuehn} and Jean-Benoit Bost \cite{\bost}, 
\cite{\bostumd}, independently. 
The computation of such arithmetic invariants via an arithmetic Lefschetz formula is
discussed  in \cite{\maillotroessler}. 
The identity (6.8) is the first arithmetic case in which the critical point $s_0$  
for $\Cal E(\tau,s)$ is not zero and 
and the leading term $\Cal E(\tau,\frac12)$ does not vanish. 
It is also the first case
in which a truely global quantity, the pairing $\langle \,\ZH(t,v),\hbom\,\rangle$ for a
horizontal  cycle $\Cal Z(t)$, must be computed; it is determined as the Faltings height 
of a CM elliptic curve, \cite{\kryII}.

{\bf3. $0$--cycles on arithmetic surfaces.} In the case $n=1$, $r=2$, 
$\wphi_2(\tau)$ is 
a generating function for $0$--cycles on the arithmetic surface $\Cal M$. 
The combination of
\cite{\annals},  joint work with Rapoport \cite{\krinvent}, and 
current joint work with Rapoport and Yang, \cite{\kryIII}, comes very close to proving the
identity
$$\degh\!(\wphi_{2}(\tau)) \qeq \Cal E'_{2}(\tau,0),\tag6.9$$
again up to an ambiguity in the constant term of the generating function due to the 
lack of a formula for $\langle\hbom,\hbom\rangle$. 
For  
$T>0$ and $p$-regular, as defined in \cite{\bourbaki}, the cycle
$\Cal Z(T)$ is a $0$-cycle concentrated in a single fiber $\Cal M_p$ for a prime $p$
determined by $T$. In this case, the computation of $\degh\!(\Cal Z(T))$ 
amounts to a counting
problem and a problem in the deformation theory of
$p$-divisible  groups. The latter is a special case of a 
deformation problem solved by Gross and Keating \cite{\grosskeating}. On
the  analytic side, the computation of the corresponding term in the central 
derivative of the Eisenstein series amounts to the {\it same} counting problem
and the computation of the central derivative of a certain Whittaker function on 
$\roman{Sp}_2(\Q_p)$. This later computation depends on the explicit formulas due to Kitaoka
\cite{\kitaoka} for the representation densities of $T$ by unimodular quadratic forms 
of rank $4+2j$, cf. \cite{\bourbaki}, section 5, for a more detailed discussion.

{\bf4. Siegel modular varieties.} $n=3$,  \cite{\krsiegel}.  The Shimura variety $M$
attached  to a rational 
quadratic space $V$ of signature $(3,2)$ is, in general, a `twisted' version of a Siegel
$3$-fold.  The canonical model $M$ over $\Q$ can
be obtained as a moduli space  of polarized abelian varieties of dimension $16$ with an
action of  a maximal order $O_C$ in the Clifford algebra of $V$. 
A model $\Cal M$ over $\Spec(\Z[N^{-1}])$, for a suitable $N$ 
can likewise be defined as a moduli space, \cite{\krsiegel}. The possible 
generating functions and their connections with Eisenstein series 
are given in the following chart, \cite{\icm}:
\smallskip
\settabs 8\columns
\+$r=1,$&\hskip -10pt$\Cal Z(t)_\Q ={\textstyle\text{Humbert}\atop\textstyle
\lower 2pt\hbox{\text{surface}}}
$,&&\hskip
-10pt 
$\hat\phi_1(\tau) = \hat\o + ? 
+\sum_{t\ne0} \ZH(t,v)\,q^t,$&&&
$\langle \hat\phi_1(\tau),\hat\o^3\rangle\overset{?}\to{ =}\ \Cal
E'_1(\tau,\frac32)$
\cr
\+$r=2,$&\hskip -10pt$\Cal Z(t)_\Q = \roman{curve}$&&\hskip -10pt 
$\hat\phi_2(\tau) = \hat\o^2 + ? 
+\sum_{T\ne0} \ZH(T,v)\,q^T,$&&&
$\langle \hat\phi_2(\tau),\hat\o^2\rangle\overset{?}\to{ =}\ \Cal
E'_2(\tau,1)$
\cr
\+$r=3,$&\hskip -10pt$\Cal Z(T)_\Q =0$--cycle,&&\hskip -10pt 
$\hat\phi_3(\tau) = \hat\o^3 + ? 
+\sum_{T\ne0} \ZH(T,v)\,q^T,$&&&
$\langle \hat\phi_2(\tau),\hat\o\rangle\overset{?}\to{ =}\ \Cal
E'_3(\tau,\frac12)$
\cr
\+$r=4,$&\hskip -10pt$\Cal Z(T)_\Q =\emptyset,$&&\hskip -10pt 
$\hat\phi_4(\tau) = \hat\o^4+ ? 
+\sum_{T\ne0} \ZH(T,v)\,q^T,$&&&
$\degh\!\hat\phi_4(\tau)\overset{?}\to{ =}\ \Cal E'_4(\tau,0).$
\cr
The Siegel Eisenstein series $\Cal E_r(\tau,s)$ 
and, conjecturally, the generating
functions 
$\hat\phi_r(\tau)$ have weight $\frac52$ and genus $r$, and 
the last column in the chart gives the `arithmetic volume formula'
of Problem 6 in each case. Some evidence for the last of 
these identities was obtained in joint work with M. Rapoport, \cite{\krsiegel}. 
In the case of a prime $p$ of good reduction a model 
of $M$ over $\Spec(\Z_p)$ is defined in \cite{\krsiegel}, and cycles are defined 
by imposing special endomorphisms. For 
$r=4$, the main results of 
\cite{\krsiegel} give a criterion for $\Cal Z(T)$ to be a 
$0$--cycle in a fiber $\Cal M_p$
and show that, when this is the case, then
$\degh\!((\Cal Z(T),0))\,q^T = \Cal E'_{4,T}(\tau,0)$. 
The calculation of the left hand side is again based on the
result of Gross and Keating, \cite{\grosskeating}.
For $r=1$, the results of \cite{\Bints}, cf., in particular, sections 5 and 6, 
are consistent with the identity in the first row, which 
involve arithmetic volumes of divisors. 

{\bf5. Divisors.} For any $n$, when $r=1$, the arithmetic volume formula predicts that 
the second term $\Cal E'_1(\tau,\frac{n}2)$ in the Laurent expansion 
of an elliptic modular 
Eisenstein series of weight $\frac{n}2+1$ at the point $s_0=\frac{n}2$
has Fourier coefficients involving the arithmetic volumes $\langle \ZH(t,v),\hbom^n\rangle$
of divisors on the integral model $\Cal M$ of $M$. The first term  
$\Cal E_1(\tau,\frac{n}2)$ in the 
Laurent expansion at this point has Fourier coefficients involving the usual
volumes of the corresponding geometric cycles.  For example, 
for a suitable choice of $V$, $\Cal E(\tau,\frac{n}2)$ 
is a familiar classical Eisenstein series, e.g., $E_2(\tau)$ (non-holomorphic),
$E_4(\tau)$, $E_6(\tau)$, etc., for $\dim(V)$ even, and Cohen's Eisenstein series 
$E_{\frac{n}2+1}(\tau)$, \cite{\cohen}, for $\dim(V)$ 
odd\footnote{Here the subscript denotes the weight rather than
the genus}.  This means that the second term
in the Laurent  expansion of such classical 
Eisenstein series should contain information 
from arithmetic 
geometry! Again, related results are obtained in \cite{\Bints}.

{\bf An Important Construction.}  
We conclude this section with an important identity which 
relates the generating function for height pairings with that 
for arithmetic degrees. 
Suppose that $n=2r-1$ is odd. Then the various conjectural 
identities above, in particular (5.7) and (6.6), lead to the formula:
$$\align
\langle \wphi_r(\tau_1),\wphi_r(\tau_2)\,\rangle 
&\ =\ \degh\big(\,\widehat{\phi_r}(\tau_1)\cdot 
\overline{\widehat{\phi}_r(\tau_2)}\,\big)\\
\nass
{}& \qeq \degh(\widehat{\phi}_{2r}\big(\pmatrix 
\tau_1&{}\\{}&-\bar\tau_2\endpmatrix\big)\,\big)\qquad\ \text{(via (5.7))}\tag6.10\\
\nass
{}& \qeq  \Cal E'_{2r}\big(\pmatrix \tau_1&{}\\{}&
-\bar\tau_2\endpmatrix,0\,\big)\qquad\qquad\,\text{(via (6.6))}
\endalign
$$ 
relating the height pairing of the series $\wphi_r(\tau)\in\CH^r(\Cal M)$ in the middle 
degree with the restriction of the central derivative of the Siegel Eisenstein series 
$\Cal E_{2r}(\tau,s)$ 
of genus $2r$ and weight $r+\frac12$. Note that this weight is always half-integral.
These series are the `incoherent' Eisenstein series discussed in \cite{\annals} and 
\cite{\bourbaki}. The conjectural identity (6.10) will be used in an essential way in the next
section. 

\vskip 1in
\centerline{\loud III}
\medskip
\medskip
\centerline{\loud Derivatives of \, L-series}

\medskip
In this section, we explain how the modularity of the arithmetic theta functions 
and the conjectural relations between their inner products and derivatives 
of Siegel Eisenstein series might be connected with higher dimensional Gross--Zagier type 
formulas expressing central derivatives of certain L-functions 
in terms of height pairings of special cycles. These formulas 
should be analogues of those connecting {\it values} of certain
L-functions to inner products of theta lifts, e.g., the Rallis inner product formula, 
\cite{\rallisinnerprod}, \cite{\li}, \cite{\krannals}, section 8, and 
Appendix II below. 

\medskip

\subheading{\Sec7. Arithmetic theta lifts}

Suppose that $f\in S^{(r)}_{\frac{n}2+1}$ is a holomorphic Siegel cusp form of weight
$\frac{n}2+1$  and genus $r$ for some subgroup $\Gamma'\subset \roman{Sp}_r(\Z)$. 
Then, assuming the existence of the generating function 
$\wphi_r$ valued in $\CH^r(\Cal M)$ and that this function is also modular for $\Gamma'$,
we can define an  {\bf arithmetic theta lift}:
$$\align
\thh_r(f) :&= \langle \, f, \wphi_r\,\rangle_{\Pet}\tag7.1\\
\nass
{}&= \int_{\Gamma'\back \H_r} f(\tau)\overline{\wphi_r(\tau)}\,
\det(v)^{\frac{n}2+1}\,d\mu(\tau) \in \CH^r(\Cal M)
\endalign
$$
where $\langle\ ,\ \rangle_{\Pet}$ is the Petersson inner product. 
Thus, we get a map
$$S_{\frac{n}2+1}^{(r)} \lra \CH^r(\Cal M), \qquad f\mapsto 
\thh_r(f).\tag7.2$$
This map is an arithmetic analogue of a theta correspondence like the Shimura lift,  
 \cite{\shimurahalf}, from 
forms of weight $\frac32$ to forms of weight $2$, which can be 
defined by integration 
against a classical theta function, \cite{\niwa}. 
For example, if $f$ is a Hecke eigenform,
then $\thh_r(f)$ will also be an Hecke eigenclass.

\subheading{\Sec8. Connections with derivatives of L-functions}

Restricting to the case $n=2r-1$, where the target is the 
arithmetic Chow group $\CH^r(\Cal M)$ 
in the middle dimension, we can compute the height pairing 
of the classes $\thh_r(f)$  
using identity (6.10) above:
$$\align
\langle\,\thh_r(f_1),\thh_r(f_2)\,\rangle &\ =\ \big\langle\ 
\langle f_1,\wphi_r\rangle_{\Pet},\langle f_2,\wphi_r\rangle_{\Pet}\,\big\rangle\\
\nass
{}&\ =\ \big\langle\,f_1\tt\bar f_2, \langle 
\wphi_r(\tau_1),\wphi_r(\tau_2)\rangle\,\big\rangle_{\Pet}\tag8.1\\
\nass
{}&\qeq \big\langle\,f_1\tt\bar f_2,\Cal E'_{2r}
\big(\pmatrix \tau_1&{}\\{}&-\bar\tau_2\endpmatrix,0\,\big)
\,\big\rangle_{\Pet}\qquad\qquad\qquad\text{by (6.10)}\\
\nass
{}&\ =\ \frac{\partial\ }{\partial s}\left\{\,\big\langle\,
f_1\tt\bar f_2,\Cal E_{2r}\big(\pmatrix \tau_1&{}\\{}&-\bar\tau_2\endpmatrix,\bar s\,\big)
\,\big\rangle_{\Pet}\right\}\bigg\vert_{s=0}
\endalign
$$
Here we use the hermitian extension of the height pairing (6.2) to $\CH^1(\Cal M)_\C$
taken to be conjugate linear in the second argument. 
Aficionados of Rankin--Selberg integrals will now 
recognize the {\bf doubling integral} 
of Rallis and Piatetski-Shapiro, \cite{\psrallis}, and, 
in classical language, of B\"ocherer, \cite{\boecherer}, and Garrett, \cite{\garrettII}, 
in the last line of (8.1). If $f_1$ and $f_2$ correspond to different 
irreducible cuspidal automorphic representations, then the integral in (8.1)
vanishes. Otherwise, one has the identity, \cite{\li}, \cite{\krannals}, 
$$\big\langle\,f_1\tt\bar f_2,\Cal E_{2r}\big(\pmatrix 
\tau_1&{}\\{}&-\bar\tau_2\endpmatrix,\bar s\,\big)
\,\big\rangle_{\Pet} = \langle \s(\P_S(s))f_1,f_2\rangle_{\Pet} \,L^S(s+\frac12,\pi)\tag8.2$$
where 
$$\align
S^{(r)}_{r+\frac12}\ni f_i &\longleftrightarrow F_i= \text{ automorphic form ($i=1$, $2$) for
$H(\A)$, for $H=SO(r+1,r)$}\\ &\matrix \text{under the analogue of the
Shimura--Waldspurger correspondence}\\
\text{ between forms of weight $r+\frac12$ on $Mp_{r}$ and forms on
$SO(r+1,r)$}\endmatrix\\
\nass
\s&=\text{ the irreducible automorphic cuspidal representation attached to $f_i$}\\
\nass
\pi&=\text{ the irreducible automorphic cuspidal representation attached to $F_i$}\\
\nass
L(s,\pi)&=\text{the degree $2r$ Langlands L-function attached to $\pi$,}\\
{}&\text{ and the standard representation 
of the L-group $H^\vee = Sp_{r}(\C)$}\\
\nass
\s(\P_S(s))&=\text{ a local `twisting' operator which gives }\\
&\qquad\text{the contribution of bad local zeta integrals.}
\endalign
$$ 
Considerations of local theta dichotomy,
\cite{\harrisksweet}, \cite{\hirschberg}, control the local root numbers 
of $L(s,\pi)$ so that the global root number is $-1$ and 
$L(\frac12,\pi)=0$.  Combining (8.1) and (8.2), we obtain the {\bf arithmetic inner
product formula}:
$$\langle\,\thh_r(f),\thh_r(f)\,\rangle \qeq B\cdot \langle f,f\rangle_{\Pet} 
\,L'(\frac12,\pi)\tag8.3$$
for a constant $B$ coming the local zeta integrals at bad primes. 
Of course, this is only conjectural! For a general discussion of what one 
expects of such central critical values, cf. \cite{\grossmot}. 

\vfill\eject

{\bf Examples.} 

{\bf 1. The Gross-Kohnen-Zagier formula.} In the case $n=r=1$, we have
$$\align
M&= \text{Shimura curve}\\
\nass
\Cal M&=\text{integral model}\\
\nass
f &= \text{weight $\frac32$}\tag8.4\\
\nass
F &= \text{corresponding form of weight $2$ (assumed a normalized newform)}\\
\nass
\pi&=\text{associated automorphic representation of $\roman{PGL}_2(\A)$}\\
\nass
\thh_1(f)&\in \CH^1(\Cal M)\\
\nass
L(s,\pi) &= L(s+\frac12,F)\\
\nass
L(s,F)&=\text{ the standard Hecke L-function (with $s\mapsto 2-s$ functional eq.).}\\
\endalign
$$
In this case, identity (8.3) becomes
$$\langle\,\thh_1(f),\thh_1(f)\,\rangle = B\cdot ||f||^2
\,L'(1,F)\tag8.5$$
This is essentially the Gross-Kohnen-Zagier 
formula, \cite{\grosskohnenzagier}, Theorem C.

{\bf 2. Curves on Siegel 3-folds.} The next example is
$n=3$ and $r=2$. Then:
$$
\align
M&=\text{ Siegel 3-fold}\\
\nass
\Cal M&=\text{ arithmetic 4-fold}\\
\nass
f&=\text{ a Siegel cusp form of weight $\frac52$ and genus $2$}\tag8.6\\
\nass
\pi&=\text{ corresponding automorphic representation of $O(3,2)$}\\
\nass
L(s,\pi)&=\text{ the degree 4 L-function of $\pi$}
\endalign
$$
The cycles $\Cal Z(T)$ in the generating function $\wphi_2(\tau)$ are now Shimura curves on
the generic fiber $M=\Cal M_\Q$,  extended to arithmetic surfaces in the arithmetic
$4$-fold $\Cal M$, and 
$$
\thh_2(f)\in \CH^2(\Cal M).\tag8.7
$$
Then identity (8.3) says that the the central 
derivative of the degree 4 L-function $L(s,\pi)$ is expressible in terms of the 
height pairing $\langle\, \thh_2(f),\thh_2(f)\rangle$ of the class $\thh_2(f)$ 
i.e., made out of the $f$--eigencomponents of
 `curves on a Siegel 3-fold'. 
Of course, the proof of such a formula by the method outlined here 
requires that we prove the relevant versions of (5.7), (6.6) and (6.10), 
and, above all, the modularity of the codimension $2$ generating function 
$\wphi_2(\tau)$. Needless to say, this remains very speculative!

{\bf 3. The central derivative of the triple product L-function.}
This case involves a slight variant of the previous pattern. 
If we take $V$ of signature $(2,2)$, we have
$$M = \cases M_1\times M_1, \quad M_1 = \text{modular curve or Shimura curve,}&\\
\nass
\text{Hilbert--Blumenthal surface}&{}\\
\nass
\text{compact Hilbert--Blumenthal type surface}&
\endcases\tag8.8
$$
where, in the two cases in the first line, the discriminant of $V$ is a square and 
$\roman{witt}(V)=2$ or $0$, respectively, while, in the second two cases, $\kay =
\Q(\sqrt{\roman{discr}(V)})$ is a real quadratic field and $\roman{witt}(V) = 1$ 
or $0$ respectively. Then, $\Cal M$ is an arithmetic $3$-fold, and, 
conjecturally, the generating function
$$\wphi_1(\tau)\in \CH^1(\Cal M)\tag8.9$$
is a modular form of weight $2$ for a subgroup $\Gamma'\subset \SL_2(\Z)$. 
Note that, on the generic fiber, the cycles $\Cal Z(t)_\Q$ are the Hirzebruch--Zagier 
curves, \cite{\krHB}.  For a cusp form $f\in S_2(\Gamma')$, we obtain a class
$$\thh_1(f)\in \CH^1(\Cal M).\tag8.10$$
Consider the trilinear form on $\CH^1(\Cal M)$ defined by
$$\langle\ \hat z_1,\hat z_2, \hat z_3\ \rangle := 
\degh\!(\hat z_1\cdot\hat z_2\cdot \hat z_3).\tag8.11$$
Then, for a triple of cusp forms of weight $2$,
$$\align
\langle\ \thh(f_1),\thh(f_2),\thh(f_3)\ \rangle &\ = \ 
\langle\ f_1f_2f_3,\langle\ \wphi_1(\tau_1),\wphi_1(\tau_2),\wphi_1(\tau_3)\ \rangle\
\rangle_{\roman{Pet}}\\
\nass
{}&\qeq \langle\ f_1f_2f_3, \degh\!(\wphi_3(
\pmatrix\tau_1&{}&{}\\{}&\tau_2&{}\\{}&{}&\tau_3\endpmatrix))\ \rangle_{\roman{Pet}}\tag8.12\\
\nass
{}&\qeq\langle\ f_1f_2f_3, \Cal E'_3(
\pmatrix\tau_1&{}&{}\\{}&\tau_2&{}\\{}&{}&\tau_3\endpmatrix,0\,)\ \rangle_{\roman{Pet}}\\
\nass
{}&\ =\ \frac{\partial\ }{\partial s}
\bigg\{\ \langle\ f_1f_2f_3, \Cal E_3(
\pmatrix\tau_1&{}&{}\\{}&\tau_2&{}\\{}&{}&\tau_3\endpmatrix,\bar s\,)\ \rangle_{\roman{Pet}}
\bigg\}\bigg\vert_{s=0}.\\
\endalign
$$
If we assume that the $f_i$'s are newforms with associated 
cuspidal automorphic representations $\pi_i$, $i=1$, $2$, $3$, then
the integral in the last line is\footnote{Almost, except that one
must actually work with the similitude group $\roman{GSp}_3$.} the integral
representation of the  triple product L-function, \cite{\garrett},
\cite{\psrallistriple}, 
\cite{\grosskudla}, \cite{\bocherersptriple}, 
$$\langle\ f_1f_2f_3, \Cal E_3(
\pmatrix\tau_1&{}&{}\\{}&\tau_2&{}\\{}&{}&\tau_3\endpmatrix,\bar s\,)\ \rangle_{\roman{Pet}}
= B(s)\,L(s+\frac12,\pi_1\tt\pi_2\tt\pi_3).\tag8.13$$

Note that the results of D. Prasad on dichotomy for local trilinear forms, 
\cite{\dprasad},  
control the local root numbers and the `target' space $V$. 

Here, in addition to the modularity of the generating function $\wphi_1(\tau)$, 
we have used the conjectural identities
$$\wphi_1(\tau_1)\cdot\wphi_1(\tau_2)\cdot\wphi_1(\tau_3) \qeq
\wphi_3(\pmatrix\tau_1&{}&{}\\{}&\tau_2&{}\\{}&{}&\tau_3\endpmatrix),\tag8.14$$
analogous to (5.7) and
$$\degh\!\wphi_3(\tau) \qeq \Cal E'_3(\tau,0),\tag8.15$$
analogous to (6.6). The equality of certain 
coefficients on the two sides of (8.15) follows from the 
result of Gross and Keating, \cite{\grosskeating}, and the formulas of 
Kitaoka, \cite{\kitaokatriple}, cf. also \cite{\krHB}. 

In fact, one of the starting points of my long crusade to establish connections between 
heights and Fourier coefficients of central derivatives of 
Siegel Eisenstein series was 
an old joint project with Gross and Zagier, of
which 
\cite{\grosskudla} was a preliminary `exercise'. 
The other was my collaboration with Michael Harris on 
Jacquet's conjecture about the central value of the
triple product L-function, 
\cite{\harriskudlaII}, \cite{\harriskudlaIII}, 
based, in turn on a long collaboration with Steve Rallis on the 
Siegel--Weil formula. And, of course, the geometric picture 
which serves as an essential guide comes from joint work with John Millson. 
I would like to thank them all, together with my current collaborators
Michael Rapoport and Tonghai Yang, for their generosity with 
their ideas, advice, encouragement, support and patience.

\define\Zhat{\widehat{\Z}}

\subheading{Appendix I: Shimura curves}

In this appendix, we illustrate some of our basic constructions in the case of 
modular and Shimura 
curves. In particular, this allows us to make a direct connection with classical 
Heegner points, one of 
the main themes of the MSRI conference. 

In the case of a rational quadratic space $V$ of signature $(1,2)$, the varieties of part I 
are the classical Shimura curves. Let $B$ be an indefinite quaternion algebra over $\Q$,
and let
$$\align
V&=\{\ x\in B\mid\tr(x)=0\ \}, \qquad Q(x) = \nu(x) = -x^2.\tag A.I.1\\
\noalign{Note that the associated bilinear form is $(x,y) = \tr(xy^\iota)$, where $x\mapsto x^\iota$ 
is the main involution on $B$. The action of $B^\times$ on $V$ by conjugation 
induces an isomorphism} 
\nass
B^\times &\isoarrow G=\roman{GSpin}(V).\tag A.I.2\\
\noalign{We fix an isomorphism} 
\nass
B_\R&=B\tt_\Q\R \isoarrow M_2(\R),\tag A.I.3
\endalign
$$
and obtain an identification
$$\Bbb P^1(\C) \setminus \Bbb P^1(\R) \isoarrow D, \qquad z \mapsto w(z) = 
\pmatrix z&-z^2\\1&-z\endpmatrix\ \mod
\C^\times.\tag A.I.4
$$  
Let $S$ be the set of the primes $p$ for which $B_p=B\tt_\Q\Q_p$ 
is a division algebra, and let $D(B) = \prod_{p\in S} p$. 
For a fixed maximal order $O_B$ of $B$, there is an isomorphism
$$\align
B(\A_f) &\isoarrow \bigg(\prod_{p\in S} B_p\bigg)\times M_2(\A_f^S),\tag A.I.5\\
\nass
O_B\tt_\Z\Zhat &\isoarrow \bigg(\prod_{p\in S} O_{B,p}\bigg)\times M_2(\Zhat^S).\\
\noalign{For an integer $N$ prime to $D(B)$, let $R$ be the Eichler order of discriminant $ND(B)$ with}
R\tt_\Z\Zhat &\isoarrow \bigg(\prod_{p\in S} O_{B,p}\bigg)
\times \{\ x\in M_2(\Zhat^S)\mid c\equiv 0\, (\roman{mod}\, N)\ \}.\tag A.I.6
\endalign
$$
Then, for the compact open subgroup $K = (R\tt_\Z\Zhat)^\times \subset G(\A_f)$, the quotient
$$
X_0^B(N):=M_K(\C) \simeq G(\Q)\back\bigg( \, D\times G(\A_f)/K\,\bigg) \simeq \Gamma\back D^+,\tag A.I.7
$$
where $\Gamma = G(\Q)^+\cap K = R^\times$, 
is the analogue for $B$ of the modular curve $X_0(N)$. Of course, when $B=M_2(\Q)$, we need to add the cusps. 
The  $0$--cycles $Z(t,\ph;K)$ are weighted combinations of CM--points. These can be described as follows. 
If we identify $V(\Q)$ with a subset of $B_\R=M_2(\R)$, then, for
$$\align
x &=\pmatrix b&2c\\-2a&-b\endpmatrix\  \in V(\Q) \subset M_2(\R),\qquad Q(x) = -(b^2-4ac),\tag A.I.8\\
\nass
D_x &= \{\ z\in \Bbb P^1(\C) \setminus \Bbb P^1(\R)\ \mid \ (x,w(z)) = -2(az^2+bz+c)=0\ \}.
\endalign
$$
Note that for general $B$, the coordinates $a$, $b$ and $c$ of $x$ need not lie in $\Q$. \hfb
For $d>0$, let
$$\O_d=\{\ x\in V\mid Q(x)=d\ \}\tag A.I.9$$
and note that, if $x_0\in \O_d(\Q)$, then 
$$\O_d(\A_f) = G(\A_f)\cdot x_0 = K\cdot \O_d(\Q).\tag A.I.10$$
By (iii) of Lemma~2.2 of \cite{\duke}, if $x\in V(\Q)$, $g\in G(\A_f)$ and $\gamma\in G(\Q)$, then, 
for the cycle defined by (2.8) above, 
$$Z(\gamma x,\gamma g;K) = Z(x,g;K).\tag A.I.11$$
Thus, 
for any $\ph\in S(V(\A_f))^K$, the weighted $0$--cycle $Z(d,\ph;K)$ on $X_0^B(N)$ 
is given by 
$$Z(d,\ph;K) = \sum_r \ph(x_r)\,Z(x_r,1;K),\tag A.I.12$$ 
with the notation of (2.11), 
where
$$\roman{supp}(\ph)\cap \O_d(\A_f) = \coprod_r K\cdot x_r, \qquad x_r\in \O_d(\Q).\tag A.I.13$$
 For example, 
if $L^\vee$ is the dual lattice of $L := R\cap V(\Q)$, there is a Schwartz function
$$\ph_\mu = \roman{char}(\mu+\widehat{L})\in S(V(\A_f))\tag A.I.14$$ 
for each coset $\mu+L$, for $\mu\in L^\vee$. Here 
$\widehat{L}=L\tt_\Z\Zhat$. The group $\Gamma$ acts on $L^\vee/L$, and each $\Gamma$--orbit 
$\Cal O$ defines a $K$--invariant weight function
$$\ph_{\Cal O} = \sum_{\mu\in \Cal O}\ph_\mu\quad \in S(V(\A_f))^K.\tag A.I.15$$
Then, 
$$Z(d,\ph_{\Cal O};K) = \sum_{\matrix \scr x\in L^\vee \cap \O_d(\Q)\\ \scr x+L\in \Cal O\\
\scr\mod \Gamma\endmatrix} \pr(D_x^+),\tag A.I.16$$
where $D_x^+ = D_x \cap D^+$ and $\pr:D^+\rightarrow \Gamma\back D^+ = X_0^B(N)$ is the projection. 
Here each point $\pr(D_x^+)$ is to be counted with multiplicity $e_x^{-1}$, where $2e_x$ is the order of the 
stablizer of $x$ in $\Gamma$. 

The Heegner cycles studied by Gross, Kohnen, and Zagier, \cite{\grosskohnenzagier} can be recovered
from this formalism in the case $B=M_2(\Q)$. Of course, we take the standard identification $B_\R=M_2(\R)$ 
and the maximal order $O_B=M_2(\Z)$. For $x\in V$, we let 
$$y =\frac12J^{-1}x= \frac12 \pmatrix {}&-1\\1&{}\endpmatrix\pmatrix b&2c\\-2a&-b\endpmatrix 
=\frac12 \pmatrix 2a&b\\b&2c\endpmatrix
=\pmatrix a&b/2\\b/2&c\endpmatrix.\tag A.I.17$$
This is the matrix for the quadratic form denoted by $[a,b,c]$ in 
\cite{\grosskohnenzagier}, p.504. Moreover, if $g\in SL_2(\Z)$, then the action of $g$ on 
$[a,b,c]$ is given by $y\mapsto {}^tg y g$, 
and this amounts to
$$x\mapsto g^{-1} x g\tag A.I.18$$
on the original $x$.
Let
$$L=\{ x=\pmatrix b&c\\-a&-b\endpmatrix \in M_2(\Z)\mid a\equiv b\equiv 0\!\!\mod 2N, 
\text{ and } c\equiv 0\!\!\mod 2\},\tag A.I.19$$
and, for a coset $r \in \Z/2N\Z$, let 
$\ph_{N,r}\in S(V(\A_f))$ be the characteristic function of the set
$$\pmatrix r&{}\\{}&-r\endpmatrix +\hat{L} \quad\subset V(\A_f).\tag A.I.20$$
Note that the function $\ph_{N,r}$ is $K$-invariant.

The set of $x$'s which contribute to $Z(d,\ph;K)$, i.e., the set 
$\text{supp}(\ph)\cap \O_d(\Q)$, is
$$\{\,x = \pmatrix b&2c\\-2a&b\endpmatrix\mid b^2-4ac =-d,\ a\equiv 0\!\!\mod N,\ 
b\equiv r\!\!\mod 2N\, \}.\tag A.I.21$$
This set is mapped bijectively 
to the set
$$\Cal Q_{N,r,d}:=\{\, y=\pmatrix a&b/2\\b/2&c\endpmatrix\mid 
b^2-4ac=-d,\ a \equiv 0\!\!\mod N 
\text{ and } b\equiv r\!\!\mod 2N\,\}\tag A.I.22$$
under the map $x\mapsto y$ of (A.I.17).
Therefore $Z(d,\ph;K)$ is precisely the image in $\Gamma_0(N)\back D^+$
of the set of roots $z$ in $D^+$, identified with the upper half plane, 
of the quadratic equations $az^2+bz+c=0$, with $[a,b,c]\in \Cal Q_{N,r,d}$.
Note that, if
$$\Cal Q_{N,r,d}^+ =\{ y\in \Cal Q_{N,r,d}\mid a>0\},\tag A.I.23$$
then
$$\Cal Q_{N,r,d}\simeq \Cal Q_{N,r,d}^+\cup \Cal Q_{N,-r,d}^+.\tag A.I.24$$
The set of roots, counted with multiplicity, 
for $[a,b,c]\in \Cal Q_{N,r,d}^+$ is denoted by $\Cal P_{-d,r}$ 
in \cite{\grosskohnenzagier}, p.542, and
$$\Cal P_{-d,r}^*=\Cal P_{-d,r}\cup \Cal P_{-d,-r},\tag A.I.25$$
where points are counted with the sum of their multiplicities in the 
two sets.
We conclude that
\proclaim{Proposition A.I.1} Fix $N$ and $r \!\!\mod 2N$, and let $\ph_{N,r}$ 
be as above. Then, for $K=K_0(N)$, as above, 
$$Z(d,\ph_{N,r};K) = \cases \Cal P_{-d,r} + \Cal P_{-d,-r}&\text{ if $-d\equiv r^2\!\!\mod 4N$}\\
0&\text{ otherwise,}
\endcases$$
as 0--cycles on $X_0(N)$.
\endproclaim

\subheading{Appendix II: Theta forms}

This Appendix contains a few more details about the theta functions valued 
in differential forms (theta forms) used to prove Theorem~3.1, and about a 
formula for cup products which is analogous to the
conjectural arithmetic inner product formula (8.3). Good references for the 
general construction of theta functions and 
automorphic forms via the Weil representation include \cite{\shintani}, \cite{\howeps},
\cite{\li}. Let
$G' =\roman{Sp}(r)$  be the symplectic group of rank $r$ over $\Q$, and let $G'_\A$ (resp. $G'_\R$) be the 
metaplectic cover\footnote{The covers are only needed when $n$ and hence $\dim(V)$ 
is odd. When $n$ is even, one can simply work with the linear groups 
$G'(\Bbbf A) = \roman{Sp}_r(\Bbbf A)$, $G'(\Bbbf R) =\roman{Sp}_r(\Bbbf R)$, etc. } of $G'(\Bbb A)$ (resp.
$G'(\Bbb R)$). Let $G'_\Q = G'(\Q) = \roman{Sp}_r(\Q)$, identified with a 
subgroup of $G'_\A$ via the unique homomorphism $\roman{Sp}_r(\Q)=G'(\Q)\rightarrow G'_\A$ 
lifting the inclusion $G'(\Q)\hookrightarrow G'(\A)$. Let $\psi$ be the character of
$\A$ which is trivial on
$\Q\cdot\widehat{\Z}$  and with $\psi_\infty(x) = e(x) = e^{2\pi i x}$. 
Then $G'_\A$ acts on the space $S(V(\A)^r)$ via the Weil representation $\o = \o_\psi$ 
determined by $\psi$, and this action 
commutes with the natural linear action of $H(\A)$, where $H=O(V)$. 
By Poisson summation, the theta distribution
$\Theta$ on
$S(V(\A)^r)$ defined by 
$$\ph\mapsto\Theta(\ph) = \sum_{x\in V(\Q)^r}\ph(x)\tag A.II.1$$ 
is invariant under the action of $\gamma\in G'_\Q$, i.e., 
$$\Theta(\o(\gamma)\ph) = \Theta(\ph).\tag A.II.2$$ 
The function
$$\theta(g',h;\ph) := \sum_{x\in V(\Q)^r} \o(g')\ph(h^{-1}x)\tag A.II.3$$
on $G'_\Q\back G'_\A \times H(\Q)\back H(\A)$ is slowly increasing. 

Returning to the Schwartz form $\ph^r_\infty(\tau)$ of (3.2), we have the basic 
identity
$$\o_\infty(g')\ph^r_\infty(\tau) = j(g',\tau)^{-(n+2)}\,\ph^r_\infty(g'(\tau)),\tag A.II.4$$
where $g'\in G'_\R$ acts on $S(V(\R)^r)$ via the Weil representation determined by 
$\psi_\infty$, and $j(g',\tau)$ is an automorphy factor with
$$j(g',\tau)^2 = \det(c\tau+d), \qquad \text{where}\quad \pmatrix a&b\\c&d\endpmatrix \tag A.II.5$$
is the image of $g'$ in $G'(\R) = \roman{Sp}_r(\R)$.
For $\ph\in S(V(\A_f)^r)$, the theta function $\theta_r(\tau,\ph)$ of (3.4) 
is given by  
$$\theta_r(\tau,\ph) = \Theta(\ph^r_\infty(\tau)\tt\ph),\tag A.II.6$$ 
up to a translation by $h\in G(\A_f)$.
Its transformation law is determined as follows.
For $\gamma\in G'_\Q$, write $\gamma = \gamma_\infty\gamma_f$, with 
$\gamma_\infty\in G'_\R$ and $\gamma_f\in G'_{\A_f}$. Then, using (A.II.4), t
$$\align
\theta_r(\gamma(\tau),\ph) & = j(\gamma_\infty,\tau)^{(n+2)}
\Theta(\o_\infty(\gamma_\infty)\ph^r_\infty(\tau)\tt\ph)\tag A.II.7\\
\nass
{}&=j(\gamma_\infty,\tau)^{(n+2)}\,\Theta(\o(\gamma)\big(\, \ph^r_\infty(\tau)\tt\o_f(\gamma_f)^{-1}\ph\big))\\
\nass
{}&=j(\gamma_\infty,\tau)^{(n+2)}\,\Theta(\ph^r_\infty(\tau)\tt\o_f(\gamma_f)^{-1}\ph)\\
\nass
{}&=j(\gamma_\infty,\tau)^{(n+2)}\,\theta_r(\tau,\o_f(\gamma_f)^{-1}\ph).
\endalign
$$
This is valid for any $\ph\in S(V(\A_f)^r)$ and any $\gamma\in G'_\Q$. To obtain a 
more traditional transformation law, 
let $K'\subset G'_{\A_f}$ be the inverse image of $\roman{Sp}_r(\widehat{\Z})$,
and let
$\Gamma'$ be the inverse image of $\roman{Sp}_r(\Z)$ in $G'_\R$. Note that, for any $\gamma'\in \Gamma'$, 
there is a unique element $k'\in K'$ such that $\gamma'k'=\gamma\in G'_\Q$. Suppose that 
$L\subset V$ is a lattice on which the quadratic form $Q$ is $\Z$--valued, and let $L^\vee$ 
be the dual lattice. Let $\widehat{L} = L\tt_\Z\widehat{\Z}$ and $\widehat{L}^\vee = L^\vee\tt_\Z\widehat{\Z}$.
Let $S_L\subset S(V(\A_f)^r)$ be the subspace of functions $\ph$ such that 
$$\roman{supp}(\ph)\subset (\widehat{L}^\vee)^r\tag A.II.8$$
and $\ph$ is constant on cosets of $\widehat{L}^r$. Then $S_L$ is preserved under the Weil 
representation 
action $\o_f$ of $K'$. Define a representation $\rho_L$ of $\Gamma'$ on $S_L$ by 
$$\rho_L(\gamma') = \o_f(k'),\tag A.II.9$$
and let $\rho_L^\vee$ be the associated representation on the dual space $S_L^\vee$.
The map 
$$\theta_r(\tau):\ \ph\mapsto \theta_r(\tau,\ph),\tag A.II.10$$
defines an element of $S_L^\vee$, and the transformation law (A.II.7) amounts to the 
classical transformation law
$$\theta_r(\gamma(\tau)) = j(\gamma',\tau)^{n+2}\,\rho_L^\vee(\gamma')\theta_r(\tau), \qquad \gamma\in
\roman{Sp}_r(\Z)\tag A.II.11$$ 
for the vector valued form $\theta_r(\tau)$, of type $(\rho_L^\vee,S_L^\vee)$ in the style of Borcherds,
\cite{\borchinventII}, 
\cite{\Bints}.

Finally, we mention the analogue in the geometric situation of the conjectural arithmetic 
inner product formula (8.3). This formula is a geometric version 
of the Rallis inner product formula, \cite{\rallisinnerprod}, 
\cite{\krannals}. A more general version of such a 
formula was used by Jian-Shu Li, \cite{\li}, to obtain nonvanishing results for the 
cohomology of locally symmetric spaces attached to classical groups. 
Here we only consider the inner product on the middle dimensional cohomology.

Suppose that $n=2r$ and that the
Shimura variety
$M_K$ is compact.  To obtain a large supply of such examples, one can work over a totally real field, as explained 
at the end of section 1. For a Siegel cusp form $f$ of genus $r$ and weight $\frac{n}2+1$, define the 
classical theta lift by 
$$\align
\theta_r(f,\ph) :&= \langle \, f, \theta_r(\cdot,\ph)\,\rangle_{\Pet}\tag A.II.12\\
\nass
{}&= \int_{\Gamma'\back \H_r} f(\tau)\,\overline{\theta_r(\tau,\ph)}\,
\det(v)^{\frac{n}2+1}\,d\mu(\tau).
\endalign
$$ 
Here $\langle\ ,\ \rangle_{\Pet}$ denotes the Petersson inner product. 
Thus, we get a map, 
$$S_{\frac{n}2+1}^{(r)} \lra \roman H^{2r}(M), \qquad f\mapsto 
[\,\theta_r(f,\ph)\,],\tag A.II.13$$
whose image lies in the subspace spanned by the cohomology 
classes of the cycles $Z(T,\ph)$. If $\ph\in S(V(\A_f)^r)^K$ is $K$--invariant, then
$[\,\theta_r(f,\ph)\,]\in \roman H^{2r}(M_K)$.
Consider the pairing of two such classes:
$$\align
(\,[\,\theta_r(f_1,\ph_1)\,],[\,\theta_r(f_2,\ph_2)\,]\,) 
&=\int_{M_K} \theta_r(f_1,\ph_1)\wedge\overline{\theta_r(f_2,\ph_2)}\\
\nass
&= \langle\ f_1\tt \bar f_2, 
\int_{M_K} \theta_r(\tau_1,\ph_1)\wedge\overline{\theta_r(\tau_2,\ph_2)}\
\rangle_{\Pet}\tag A.II.14\\
\nass
&= \langle\ f_1\tt\bar f_2, \int_{M_K} \theta_{n}(
\pmatrix\tau_1&{}\\{}&-\bar\tau_2\endpmatrix,\ph_1\tt\bar\ph_2)\
\rangle_{\Pet}\\
\nass
{}&= \vol(M_K,\O^n)\cdot\langle\ f_1\tt\bar f_2, E_{2n}(
\pmatrix\tau_1&{}\\{}&-\bar\tau_2\endpmatrix,\frac12,\ph_1\tt\bar\ph_2)\
\rangle_{\Pet}
\endalign
$$
by (3.7) and (4.4). In contrast to the situation in (8.1) and 
(8.2), the doubling integral occuring in the last line here
involves Siegel modular cusp forms  $f_i\in S_{r+1}^{(r)}$ of {\it integral} weight. The integral is 
zero unless the $f_i$'s correspond to the same irreducible cuspidal automorphic 
representation $\pi$ of $\roman{Sp}_r(\A)$, in which case the doubling integral 
gives, \cite{\krannals}, (7.2.8), p.69,
$$\langle\ f_1\tt\bar f_2, E_{2n}(
\pmatrix\tau_1&{}\\{}&-\bar\tau_2\endpmatrix,\bar s,\ph_1\tt\bar\ph_2)\
\rangle_{\Pet}
= \frac{1}{b_{2r}^S(s,\chi_V)}\cdot \langle \pi(\P_S(s))f_1,f_2\rangle_{\Pet} \,L^S(s+\frac12,\pi)
\tag A.II.15
$$
where 
$$\align
L(s,\pi)&=\text{the degree $2r+1$ Langlands L-function attached to $\pi$,}\tag A.II.16\\
{}&\text{ and the standard representation 
of the L-group $\roman{Sp}_{r}^\vee = SO(2r+1,\C)$}\\
\nass
\pi_S(\P(s))&=\text{a convolution operator, determined by $\ph_1\tt\bar\ph_2$, which gives}\\
{}&\text{the contribution of bad local zeta integrals.}\\
\nass
b_{2r}(s,\chi_V) &= L(s+r+\frac12,\chi_V)\cdot\prod_{k=1}^r \zeta(2s+2k-1),
\endalign
$$ 
and the superscript $S$ indicates that the Euler factors for places in a 
finite set $S$ (including the archimedean places) have been omitted. 
We refer the reader to sections 7 and 8 of \cite{\krannals} for more details. 
Note that, in the arithmetic situation of section 8, 
the factor analogous to $b_{2r}(s,\chi_V)$ is included in the
normalization of the  Eisenstein series $\Cal E_{2r}(\tau,s)$ and hence does not show up in (8.2). 
Thus, we obtain the analogue of (8.3) for the cup product:
$$(\,[\,\theta_r(f,\ph)\,],[\,\theta_r(f,\ph)\,]\,)  
=B\cdot\langle f,f\rangle_{\Pet} \,L(1,\pi),\tag A.II.17$$
where we have lumped various constant factors in $B$. 

We also note that, in our anisotropic case, Theorem~4.1 gives
$$\vol(M_K,\O^n) \cdot E_n(\tau,\frac12,\ph) = \sum_{T\ge 0} \vol(Z(T,\ph))\,q^T\tag A.II.18$$
and, for $r_1+r_2=n$, the analogue of the conjectural identity (6.10) is the formula
$$([\theta_{r_1}(\tau_1,\ph_1)], [\theta_{r_2}(\tau_2,\ph_2)]) = 
\vol(M_K,\O^n) \cdot E_n(\pmatrix \tau_1&{}\\{}&-\bar\tau_2\endpmatrix,\frac12,\ph_1\tt\bar\ph_2),\tag A.II.19$$
expressing the cup product of the cohomology generating functions 
valued in $\roman{H}^{2r_1}(M_K)$ and $\roman{H}^{2r_2}(M_K)$ as 
the pullback to $\H_{r_1}\times\H_{r_2}$ of this Eisenstein series. 

\vskip.5in


\redefine\vol{\oldvol}

\bye